\input amstex
\input Amstex-document.sty

\pageno 615

\def\shorttitle{Branching Problems of Unitary Representations}
\def\shortauthor{T. Kobayashi}

\topmatter
\title\nofrills{\boldHuge Branching Problems of Unitary Representations}
\endtitle

\author \Large Toshiyuki Kobayashi* \endauthor

\thanks *RIMS, Kyoto University, Kyoto 606-8502, Japan.  E-mail: toshi\@kurims.kyoto-u.ac.jp \endthanks

\abstract\nofrills \centerline{\boldnormal Abstract}

\vskip 4.5mm

The irreducible decomposition
 of a unitary representation often contains continuous spectrum
 when restricted to a non-compact subgroup.
The author singles out a nice class of branching problems
 where each irreducible summand occurs discretely with finite multiplicity
  (admissible restrictions).
Basic theory and new perspectives of admissible restrictions
 are presented from both analytic and algebraic view points.
We also discuss some applications of admissible restrictions
 to modular varieties and $L^p$-harmonic analysis.

\vskip 4.5mm


\endabstract
\endtopmatter

\document

\baselineskip 4.5mm \parindent 8mm

\specialhead \noindent \boldLARGE 1. Introduction \endspecialhead

Let $\pi$ be an irreducible unitary representation of a group $G$.
A {\bf{branching law}} is the irreducible decomposition of
  $\pi$ when restricted to a subgroup $G'$:
$$
     \pi|_{G'}
     \simeq \int_{\widehat{G'}}^{\oplus} m_{\pi}(\tau)\tau
     \ d \mu (\tau)
     \qquad \text{(a direct integral)}.
\tag 1.1
$$
Such a decomposition is unique, for example, if $G'$ is a
  reductive Lie group,
and the {\bf{multiplicity}}
 $m_{\pi}:\widehat {G'} \to {\Bbb N} \cup\{\infty\}$
 makes sense as a measurable function
  on the unitary dual $\widehat {G'}$.

Special cases of {\bf branching problems} include
 (or reduce to) the followings:
Clebsch-Gordan coefficients, Littlewood-Richardson rules,
decomposition of tensor product representations,
character formulas, Blattner formulas,
Plancherel theorems for homogeneous spaces, description
of breaking symmetries in quantum mechanics,
theta-lifting in automorphic forms, etc.
The restriction of unitary representations serves
 also as a method to study discontinuous groups
 for non-Riemannian homogeneous spaces
 (e.g. \cite{Mg, Oh}).

Our interest is in the branching problems
 for (non-compact) reductive Lie groups $G \supset G'$.
In this generality,
 there is no known algorithm to find branching laws.
Even worse,
 branching laws usually contain both
 discrete and continuous spectrum
 with possibly infinite multiplicities
 (the multiplicity is infinite, for example, in the decomposition
 of the tensor product of two principal series representations
 of $SL(n, \Bbb {C})$ for $n \ge 3$,
  \cite{Ge-Gr}).

The author introduced the notion
 of {\bf admissible restrictions}
 and  {\bf infinitesimal discrete decomposability}
 in \cite{$\text{Ko}_5$} and  \cite{$\text{Ko}_9$}, respectively,
 seeking for a good framework of branching problems,
 in which we could
 expect especially a simple and detailed study of branching laws,
 which in turn might become powerful methods in other fields as well
  where restrictions of representations naturally arise.

The criterion in Theorem~B indicates that
 there is a fairly rich examples of admissible restrictions;
 some are known and the others are new.
In this framework, a number of explicit branching laws
 have been newly found
 (e.g\. \cite{D-Vs, $\text{Gr-W}_{1,2}$, Hu-P-S, $\text{Ko}_{1,3,4,8}$,
 $\text{Ko-\O}_{1,2}$, $\text{Li}_2$, $\text{Lo}_{1,2}$, X}).
The point here is that branching problems become accessible
 by algebraic techniques if there is no continuous spectrum.

The first half of this article surveys briefly
 a general theory of admissible restrictions
 both from analytic and algebraic view points (\S 2, \S 3).
For the simplicity of exposition,
 we restrict ourselves to unitary representations,
 although a part of the theory can be generalized to
 non-unitary representations.
The second half discusses some applications
 of discretely decomposable restrictions.
The topics range from representation theory itself (\S 4) to
 some other fields such as $L^p$-analysis on
  non-symmetric homogeneous spaces (\S 5)
 and
 topology of modular varieties (\S 6).

\specialhead \noindent \boldLARGE 2.  Admissible restrictions to
subgroups
\endspecialhead

Let $G'$ be a subgroup of $G$,
 and $\pi \in \widehat{G}$.
 In light of (1.1),
  we introduce:
\definition{Definition~2.1} {\it
We say the restriction $\pi|_{G'}$ is {\bf{$G'$-admissible}}
 if it decomposes {\bf discretely}
 and the multiplicity $m_{\pi}(\tau)$ is {\bf finite}
 for any $\tau \in \widehat {G'}$.}
\enddefinition

One can easily prove the following assertion:
\proclaim{Theorem~A \ {\rm (\cite{$\text{Ko}_5$, Theorem~1.2})}}
Let $G \supset G' \supset G''$ be a chain of groups,
 and $\pi \in \widehat{G}$.
If the restriction $\pi|_{G''}$ is $G''$-admissible,
 then $\pi|_{G'}$ is $G'$-admissible.
\endproclaim

Throughout this article,
 we shall treat the setting as below:
\proclaim{Definition~2.2} {\it We say $(G, G')$ is a {\bf pair of
reductive Lie groups} if}

1) {\it $G$ is a real reductive linear Lie group or its finite
cover,  and}

2) {\it $G'$ is a closed subgroup, and is reductive in $G$.}

\noindent
Then, we shall fix maximal compact subgroups
 $K \supset K'$ of $G \supset G'$, respectively.
\endproclaim

A typical example is a {\bf reductive symmetric pai} $(G, G')$,
 by which we mean that $G$ is as above and that
 $G'$ is an open subgroup of the set $G^\sigma$
   of the fixed points of an involutive automorphism $\sigma$ of $G$.
For example,
 $(G, G') = (GL(n, \Bbb C), GL(n, \Bbb R))$,
 $(SL(n, \Bbb R), SO(p,n-p))$
 are the cases.

Let $(G, G')$ be a pair of reductive Lie groups.
Here are previously known examples of {\bf admissible restrictions}:
\example{Example~2.3}
The restriction $\pi|_{G'}$ is $G'$-admissible in the following cases:
\item{1)}
 (Harish-Chandra's admissibility theorem)
 $\pi \in \widehat{G}$ is arbitrary and $G'=K$.
\item{2)}
 (Howe, \cite{$\text{Ho}_1$}) $\pi$ is the Segal-Shale-Weil representation
 of the metaplectic group $G$,
 and its subgroup $G'=G_1' G_2'$ forms a dual pair
 with $G_1'$ compact.
\endexample
 In these examples, either the subgroup $G'$ or the representation $\pi$
 is very special,
 namely, $G'$ is compact or $\pi$ has a highest weight.
Surprisingly, without such assumptions,
it can happen that the restriction $\pi|_{G'}$ is $G'$-admissible.
The following criterion
  asserts that the \lq\lq  balance\rq\rq\
 of $G'$ and $\pi$ is crucial to the $G'$-admissibility.
\proclaim{Theorem~B (criterion for admissible restrictions,
 \cite{$\text{Ko}_7$})}
Let $G \supset G'$ be a pair of reductive Lie groups,
 and $\pi \in \widehat G$.
If
$$
   \operatorname{Cone}(G') \cap \operatorname{AS}_K(\pi)=\{0\},
\tag 2.1
$$
 then the restriction $\pi|_{K'}$ is $K'$-admissible.
 In particular,
 the restriction $\pi|_{G'}$ is $G'$-admissible,
namely,
 decomposes discretely
 with finite multiplicity.
\endproclaim
 A main tool of the proof of Theorem~B is the microlocal study of characters
  by using the singularity spectrum of hyperfunctions.
 The idea goes back to Atiyah, Howe, Kashiwara and Vergne
  \cite{A, $\text{Ho}_2$, Ks-Vr} in the late '70s.
The novelty of Theorem~B is to establish a framework of
{\bf admissible restrictions} with a number of new examples of interest,
which rely on a deeper understanding of the unitary dual
 developed largely in the '80s (see \cite{Kn-Vo} and references therein).

Let us briefly explain the notation used in Theorem~B.
We write $\frak k_0' \subset \frak k_0$
  for the Lie algebras of $K' \subset K$, respectively.
Take a Cartan subalgebra $\frak t_0$ of $\frak k_0$.
Then,
 $\operatorname{AS}_K(\pi)$ is the {\bf asymptotic $K$-support} of $\pi$
   (\cite{Ks-Vr}),
  and $\operatorname{Cone}(G')$ is defined
   as
$$
 \operatorname{Cone}(G')
  := \sqrt{-1} (\frak t_0^* \cap \operatorname{Ad}^*(K)({\frak k_0'}^\perp)).
 \tag 2.2
$$
By definition, both
 $\operatorname{AS}_K(\pi)$ and $\operatorname{Cone}(G')$
  are closed cones in $\sqrt{-1} \frak t_0^*$.
\example{Example~2.4}
 If $G' = K$, then the assumption (2.1) is automatically fulfilled
 because $\operatorname{Cone}(G') = \{0\}$.
The conclusion of Theorem~B in this special case
 is nothing but {\bf Harish-Chandra's admissibility theorem}
  (Example~2.3~(1)).
\endexample
To apply Theorem~B for non-compact $G'$,
 we rewrite the assumption (2.1) more explicitly
  in specific settings.
On the part $\operatorname{Cone}(G')$, we mention:
\example{Example~2.5}
$\operatorname{Cone}(G')$ is a linear subspace
 $\sqrt{-1}{(\frak t_0^*)}^{-\sigma}$ (modulo the Weyl group)
 if $(G, G')$ is a reductive symmetric pair
 given by an involution $\sigma$.
Here,
 we have chosen a Cartan subalgebra $\frak t_0$
  to be maximally $\sigma$-split.
\endexample
On the part $\operatorname{AS}_K(\pi)$,
 let us consider a unitary representation $\pi_\lambda$
  which is \lq\lq attached to\rq\rq\
 an elliptic coadjoint orbit
  $\Cal O_\lambda :=
  \operatorname{Ad}^*(G) \lambda$,
  in the orbit philosophy due to Kirillov-Kostant.
This representation is a unitarization of
 a Zuckerman-Vogan module $A_{\frak q}(\lambda)$
 after some $\rho$-shift,
  and can be realized in
 the Dolbeault cohomology group on $\Cal O_\lambda$
 by the results of Schmid and Wong.
(Here, we adopt the same polarization and normalization
 as in a survey \cite{$\text{Ko}_4$, \S 2},
 for the {\bf geometric quantization}
  $\Cal O_\lambda \Rightarrow \pi_\lambda$.)
We note that $\pi_\lambda \in \widehat G$ for \lq\lq most\rq\rq\ $\lambda$.
Let $\frak g = \frak k + \frak p$ be the complexification of a
 Cartan decomposition of the Lie algebra $\frak g_0$ of $G$.
We set
$$
     \Delta_{\lambda}^{+}(\frak p)
     :=\{\alpha \in \Delta(\frak p, \frak t):
            \langle \lambda, \alpha\rangle >0\},
            \quad
\text{ for } \lambda \in \sqrt{-1}\frak t_0^*.
$$
The original proof (see \cite{$\text{Ko}_5$})
of the next theorem was based on
  an algebraic method without using microlocal analysis.
Theorem~B gives a simple and alternative proof.
\proclaim{Theorem~C \rm (\cite{$\text{Ko}_5$})}
Let $\pi_{\lambda}\in \widehat G$
 be attached to an elliptic coadjoint orbit $\Cal O_\lambda$.
If
$$
     \Bbb R \operatorname{-span}\Delta_{\lambda}^{+}(\frak p)
     \cap
     \operatorname{Cone}(G')=\{0\},
  \tag 2.2
$$
then the restriction $\pi_{\lambda}|_{G'}$
 is $G'$-admissible.
\endproclaim
Let us illustrate Theorem~C in Examples~2.6 and 2.7
 for non-compact $G'$.
For this,
 we note that a maximal compact
 subgroup $K$ is sometimes of the form $K_1 \times K_2$ (locally).
This is the case if $G/K$ is a Hermitian symmetric space
 (e.g. $G = Sp(n, \Bbb R), SO^*(2n), SU(p,q)$).
It is also the case if
 $G=O(p,q)$, $Sp(p,q)$, etc.
\example{Example~2.6 {\rm ($K \simeq K_1 \times K_2$)}}
Suppose $K$ is (locally) isomorphic to the direct
product group $K_1 \times K_2$.
Then,
 {\sl the restriction $\pi_{\lambda}|_{G'}$
 is $G'$-admissible
  if $\lambda|_{\frak t \cap \frak k_2}=0$ and $G' \supset K_1$.}
So does the restriction $\pi|_{G'}$
 if $\pi$ is any subquotient of a coherent continuation
 of $\pi_{\lambda}$.
This case was a prototype of $G'$-admissible restrictions $\pi|_{G'}$
   (where $G'$ is non-compact and  $\pi$ is a non-highest weight module)
 proved in 1989 by the author
 (\cite{$\text{Ko}_1$; $\text{Ko}_2$, Proposition~4.1.3}),
 and was later generalized to Theorems~B and C.
Special cases include:
\item{(1)}\enspace
$K_1 \simeq \Bbb T$,
 then $\pi$ is a unitary highest weight module.
The admissibility of the restrictions $\pi|_{G'}$ in this case
 had been already known in '70s
 (see Martens \cite{Mt}, Jakobsen-Vergne \cite{J-Vr}).
\item{(2)}\enspace
$K_1 \simeq SU(2)$,
 then $\pi_{\lambda}$ is
 a quaternionic discrete series.
Admissible restrictions $\pi|_{G'}$ in this case
 are especially studied by Gross and Wallach
 \cite{$\text{Gr-W}_1$} in '90s.
\item{(3)}\enspace
$K_1 \simeq O(q), U(q), Sp(q)$.
Explicit branching laws of the restriction $\pi_\lambda|_{G'}$
  for singular $\lambda$ are given in \cite{$\text{Ko}_3$, Part I}
 with respect to the vertical inclusions of the diagram below
 (see also \cite{$\text{Ko}_1, \text{Ko}_5$}
 for those to horizontal inclusions).
%
$$
  \alignat 5
  &\hphantom{MM}O(4p,4q)  &&\supset
  &&\hphantom{MM}U(2p,2q)  &&\supset
  &&\hphantom{MM}Sp(p,q)
\\
  &\hphantom{MMm}\cup      &&
  &&\hphantom{MMii}\cup    &&
  &&\hphantom{MMm}\cup
\\
  &O(4r) \times O(4p-4r,4q) &&\supset\hphantom{i} &&U(2r)\times U(2p-2r,2q)   &&\supset\hphantom{i}
&&Sp(r) \times Sp(p-r,q)
  \endalignat
$$
 \endexample
\example{Example~2.7 {\rm (conformal group)}}
There are 18 series of irreducible unitary representations of $G :=U(2,2)$
  with regular integral infinitesimal characters.
Among them, 12 series (about \lq\lq 67\% \rq\rq\ !)
  are $G'$-admissible when restricted to
   $G':= Sp(1,1)$.
%
\endexample
The assumption in Theorem~B is in fact necessary.
By using the technique of symplectic geometry,
 the author proved the converse statement of Theorem~B:
\proclaim{Theorem~D\ {\rm (\cite{$\text{Ko}_{13}$})}}
Let $G \supset G'$ be a pair of reductive Lie groups,
 and $\pi \in \widehat G$.
If the restriction $\pi|_{K'}$ is $K'$-admissible,
 then
$\operatorname{Cone}(G') \cap \operatorname{AS}_K(\pi)=\{0\}$.
\endproclaim

\specialhead \noindent \boldLARGE 3.  Infinitesimal discrete
decomposability
\endspecialhead

The definition of admissible restrictions (Definition~2.1)
 is \lq\lq analytic\rq\rq,
  namely, based on the direct integral decomposition (1.1)
  of unitary representations.
Next, we
 consider discrete decomposable restrictions
 by a purely algebraic approach.

\definition{Definition~3.1 \ {\rm (\cite{$\text{Ko}_9$, Definition~1.1})}}
{\it Let $\frak g$ be a Lie algebra. We say a $\frak g$-module $X$
is
 {\bf {discretely decomposable}}
 if there is an increasing sequence of $\frak g$-submodules of finite length:
$$
  X = \bigcup_{m=0}^\infty X_m,
  \quad  X_0 \subset X_1 \subset X_2 \subset \dotsb.
\tag 3.1
$$
We note that $\dim X_m = \infty$ in most cases below.}
\enddefinition
Next, consider the restriction of group representations.
\definition{Definition~3.2} {\it
Let $G \supset G'$ be a pair of reductive Lie groups,
 and $\pi \in \widehat G$.
We say that the restriction $\pi|_{G'}$ is
 {\bf {infinitesimally discretely decomposable}}
 if the underlying $(\frak g, K)$-module $\pi_K$ is discretely decomposable
 as a $\frak g'$-module.}
\enddefinition
The terminology \lq\lq discretely decomposable\rq\rq\ is named
after the following fact: \proclaim{Theorem~E {\rm
(\cite{$\text{Ko}_9$})}} Let $(G,G')$ be a pair of reductive Lie
groups,
 and $\pi_K$ the underlying $(\frak g, K)$-module of $\pi \in \widehat{G}$.
Then {\rm (i)} and {\rm (ii)} are equivalent:
\item{{\rm i)}}\enspace
 The restriction $\pi|_{G'}$
 is infinitesimally discretely decomposable.
\item{{\rm ii)}}\enspace
The $(\frak g, K)$-module $\pi_K$ has a {\bf discrete branching law}
 in the sense that $\pi_K$ is isomorphic to an algebraic direct sum
 of irreducible $(\frak g', K')$-modules.
\endproclaim

Moreover,
 the following theorem holds:
\proclaim{Theorem~F
 (infinitesimal $\Rightarrow$ Hilbert space decomposition;
  \cite{$\text{Ko}_{11}$})}
Let $\pi \in \widehat{G}$.
If the restriction $\pi|_{G'}$
 is infinitesimally discretely decomposable,
 then the restriction $\pi|_{G'}$ decomposes without continuous spectrum:
$$
   \pi|_{G'} \simeq \underset{\tau \in \widehat{G'}}\to{{\sum}^\oplus}
    m_\pi(\tau) \tau
  \quad \text{(a discrete direct sum of Hilbert spaces)}.
\tag 3.2
$$
\endproclaim
At this stage,
 the multiplicity
$
 m_\pi(\tau):= \dim \operatorname{Hom}_{G'}(\tau, \pi|_{G'})
$
   can be infinite.

However,
 for a reductive symmetric pair $(G, G')$,
 it is likely that
 the multiplicity of discrete spectrum is finite
 under the following assumptions, respectively.
\item{(3.3)}\enspace
 $\pi$ is a discrete series representation for $G$.
\item{(3.4)}\enspace
The restriction $\pi|_{G'}$ is infinitesimally discretely decomposable.
\proclaim{Conjecture~3.3 {\rm (Wallach, \cite{X})}}
 $m_{\pi}(\tau)<\infty$ for any $\tau \in \widehat{G'}$
 if (3.3) holds.
\endproclaim
\proclaim{Conjecture~3.4 {\rm (\cite{$\text{Ko}_{11}$, Conjecture~C})}}
 $m_{\pi}(\tau)<\infty$ for any $\tau \in \widehat{G'}$
 if (3.4) holds.
\endproclaim
We note that Conjecture~3.4 for compact $G'$
 corresponds to Harish-Chandra's admissibility theorem.
A first affirmative result for general non-compact $G'$
 was given in \cite{$\text{Ko}_9$},
  which asserts that Conjecture~3.4 holds
  if $\pi$ is attached to an elliptic coadjoint orbit.
A special case of this assertion is:
\proclaim{Theorem~G {\rm (\cite{$\text{Ko}_{9}$})}}
$m_{\pi}(\tau)<\infty$ for any $\tau \in \widehat{G'}$
 if both (3.3) and (3.4) hold.
\endproclaim
In particular,
 Wallach's Conjecture~3.3 holds in the discretely decomposable case.
We note that
 an analogous finite-multiplicity statement fails
   if {\bf continuous spectrum}
 occurs in the restriction $\pi|_{G'}$
 for a reductive symmetric pair $(G, G')$:
\proclaim{Counter Example~3.5 {\rm(\cite{$\text{Ko}_{11}$})}}
{\rm $m_{\pi}(\tau)$ can be $\infty$
 if neither (3.3) nor (3.4) holds.}
\endproclaim
Recently, I was informed by
 Huang and Vogan that they proved Conjecture~3.4 for any $\pi$ \cite{Hu-Vo}.

A key step of Theorem~G
 is to deduce the $K'$-admissibility
 of the restriction $\pi|_{K'}$ from the discreteness assumption (3.4),
  for which we employ Theorem~H below. %
Let us explain it briefly.
We write $\Cal V_{\frak g}(\pi)$
 for the {\bf associated variety} of the underlying $(\frak g, K)$-module
  of $\pi$ (see \cite{Vo}),
 which is an algebraic variety contained in
 the {\bf nilpotent cone} of $\frak g^*$.
Let $\operatorname{pr}_{\frak g \to \frak g'}: \frak g^* \to
({\frak g'})^*$ be the projection
 corresponding to $\frak g' \subset \frak g$.
Here is a necessary condition for infinitesimal discrete
 decomposability:
\proclaim{Theorem~H (criterion for discrete decomposability
 {\rm {\cite{$\text{Ko}_9$, Corollary~3.4}}})}
Let $\pi\in \widehat{G}$.
If the restriction $\pi|_{G'}$ is infinitesimally discretely decomposable,
 then
  $\operatorname{pr}_{\frak g\to\frak g'}(\Cal V_{\frak g}(\pi))$ is contained
 in the nilpotent cone of $(\frak g')^*$.
\endproclaim

We end this section with a useful information on irreducible
 summands.
\proclaim{Theorem~I (size of irreducible summands,
{\rm{\cite{$\text{Ko}_9$}}})}
Let $\pi \in \widehat{G}$.
If the restriction $\pi|_{G'}$ is infinitesimally discretely decomposable,
 then any irreducible summand has the same
 associated variety,
   especially,  the same Gelfand-Kirillov dimension.
\endproclaim
Here is a special case of Theorem~I:
\example{Example~3.6 {\rm (highest weight modules,
{\rm{\cite{$\text{N-Oc-T}$}}})}}
Let $G$ be the metaplectic group,
 and $G'= G_1' G_2'$ is a dual pair with $G_1'$ compact.
Let $\theta(\sigma)$ be an irreducible unitary highest weight module
 of $G_2'$ obtained as
  the {\bf theta-correspondence} of $\sigma \in \widehat{G_1'}$.
Then the associated variety of $\theta(\sigma)$ does
 not depend on $\sigma$, but only on $G_1'$.
\endexample
An analogous statement to Theorem~I fails
 if there exists continuous spectrum
 in the branching law $\pi|_{G'}$
  (see \cite{$\text{Ko}_{11}$} for counter examples).

\specialhead \noindent \boldLARGE 4. Applications to
representation theory
\endspecialhead

So far,
 we have explained basic theory of discretely decomposable restrictions
 of unitary representations
 for reductive Lie groups $G \supset G'$.
Now,
 we ask what discrete decomposability
 can do for representation theory.
Let us clarify advantages of admissible restrictions,
 from which the following applications (and some more)
  have been brought out and seem to be promising furthermore.

\item{1)}\enspace
Study of $\widehat {G'}$ as irreducible summands of $\pi|_{G'}$.

\item{2)}\enspace
Study of $\widehat G$ by means of the restrictions to subgroups $G'$.

\item{3)}\enspace
Branching laws of their own right.

\subhead{4.1}\endsubhead
{}From the view point of the study of $\widehat {G'}$
 (smaller group),
 one of advantages of admissible restrictions
 is that each irreducible summand of the branching law $\pi|_{G'}$
 gives an explicit construction of an element of $\widehat{G'}$.

Historically,
 an early success of this idea (in '70s and '80s)
 was the construction of irreducible highest weight modules
(Howe, Kashiwara-Vergne, Adams, $\cdots$).
A large part of these modules can be constructed
 as irreducible summands of discrete branching laws of the Weil representation
 (see Examples~2.3~(2) and 3.6).

This idea works also
 for non-highest weight modules.
As one can observe from the criterion in Theorem~B,
 the restriction $\pi|_{G'}$ tends to be
 discretely decomposable,
 if $\operatorname{AS}_K(\pi)$ is \lq\lq{small}\rq\rq.
In particular,
 if $\pi$ is a {\bf{minimal representation}}
 in the sense that its annihilator is the Joseph ideal,
 then a result of Vogan implies that $\operatorname{AS}_K(\pi)$ is
 one dimensional.
Thus,
 there is a good possibility of finding subgroups $G'$
 such that $\pi|_{G'}$ is $G'$-admissible.
This idea was used to construct \lq\lq small\rq\rq\ representations
 of subgroups $G'$ %
  by Gross-Wallach \cite{$\text{Gr-W}_1$}.
In the same line,
 discretely decomposable branching laws for non-compact $G'$
 are used also
 in the theory of automorphic forms
 for exceptional groups by J-S. Li \cite{$\text{Li}_2$}.

\subhead
{4.2}
\endsubhead
{}From the view point of the study of $\widehat G$
 (larger group),
 one of advantages of admissible restrictions
 is to give a clue to a detailed study of representations
 of $G$ by means of discrete branching laws.

Needless to say,
 an early success in this direction
 is the theory of $(\frak g, K)$-modules
 (Lepowsky, Harish-Chandra, $\cdots$).
The theory relies heavily on Harish-Chandra's admissibility theorem
 (Example~2.3~(1)) on the restriction of $\pi$ to $K$.

Instead of a maximal compact subgroup $K$,
 this idea applied to a non-compact subgroup $G'$
 still works,
 especially in the study of \lq\lq small\rq\rq\
  representations of $G$.
In particular, this approach makes sense if
 the $K$-type structure is complicated
 but the $G'$-type structure is less complicated.
Successful examples in this direction include:
\item{1)}\enspace
To determine an explicit condition on $\lambda$
 such that a Zuckerman-Vogan module $A_{\frak q}(\lambda)$ is non-zero,
  where we concern with the parameter $\lambda$
 outside the good range.
In the setting of Example~2.6~(3),
 the author found in \cite{$\text{Ko}_2$} a combinatorial formula
 on $K_1$-types of $A_{\frak q}(\lambda)$ and determined explicitly when
  $A_{\frak q}(\lambda) \neq 0$.
The point here
 is that the computation of $K$-types of $A_{\frak q}(\lambda)$
 is too complicated to carry out because a lot of cancellation occurs
  in the generalized Blattner formula,
   while $K_1$-type formula
   (or $G'$-type formula for some non-compact subgroup $G'$)
    behaves much simpler in this case.

\item{2)}\enspace
To study a fine structure of standard representations.
For example, Lee and Loke \cite{Le-Lo}
 determined the Jordan-H\"older series and the unitarizability of
  subquotients
  of certain degenerate non-unitary principal series representations $\pi$,
   by using
  $G'$-admissible restrictions
   for some non-compact reductive subgroup $G'$.
  Their method works successfully even
 in the case where $K$-type  multiplicity of $\pi$ is not one.

\subhead
{4.3}
\endsubhead
{}From the view point
 of finding explicit branching law,
 an advantage of admissible restrictions
  is that one can employ algebraic techniques
  because of the lack of continuous spectrum.
A number of explicit branching laws are newly found
  (e.g\. \cite{D-Vs, $\text{Gr-W}_{1,2}$, Hu-P-S, $\text{Ko}_{1,3,4,8}$,
 $\text{Ko-\O}_{1,2}$, $\text{Li}_2$, $\text{Lo}_{1,2}$, X})
 in the context of admissible restrictions to non-compact reductive subgroups.
A mysterious feature is that \lq\lq different series\rq\rq\ of
 irreducible representations
 may appear in discretely decomopsable branching laws
 (see \cite{$\text{Ko}_{5}$, p\.184} for a precise meaning),
 although all of them have the same Gelfand-Kirillov dimensions (Theorem~I).

\specialhead \noindent \boldLARGE
 5. New discrete series for homogeneous spaces \endspecialhead

Let  $G \supset H$ be a pair of reductive Lie groups.
Then, there is a $G$-invariant Borel measure on the homogeneous
 space $G/H$,
 and one can define naturally a unitary representation
 of $G$ on the Hilbert space $L^2(G/H)$.
\proclaim{Definition~5.1} We say $\pi$ is a {\bf{discrete series
representation}} for $G/H$,
 if $\pi \in \widehat{G}$ is realized as a subrepresentation of
 $L^2(G/H)$.
\endproclaim

A discrete series representation corresponds
 to a discrete spectrum
 in the Plancherel formula for the homogeneous space $G/H$.
One of basic problems in non-commutative harmonic analysis is:
\proclaim{Problem 5.2} {\rm 1)\enspace
 Find a condition on the pair of groups $(G,H)$
 such that there exists a discrete series representation
 for the homogeneous space $G/H$.

{\rm 2)}\enspace If exist,
 construct discrete series representations.
\endproclaim
Even the first question has not found a final answer
 in the generality that $(G,H)$ is a pair
 of reductive Lie groups.
Here are some known cases:
\example{Example~5.3}
Flensted-Jensen, Matsuki and Oshima  proved
 in '80s that
discrete series representations for a reductive symmetric space $G/H$
 exist
 if and only if
$$
 \operatorname{rank} G/H=\operatorname{rank} K/(H \cap K).
\tag 5.1
$$
This is a generalization of Harish-Chandra's condition,
 $\operatorname{rank} G = \operatorname{rank} K$,
 for a group manifold $G \times G/\operatorname{diag}(G) \simeq G$
 (\cite{FJ, Mk-Os}).
\endexample

Our strategy to attack Problem~5.2
 for more general (non-symmetric) homogeneous spaces $G/H$
 consists of two steps:
 \item{1)}
  To embed $G/H$ into a larger homogeneous space
 $\widetilde{G}/\widetilde{H}$, on which harmonic analysis is well-understood
  (e.g. symmetric spaces).
\item{2)}
 To take functions belonging to a discrete series representation
 $\Cal H \ (\hookrightarrow L^2(\widetilde{G}/\widetilde{H}))$,
 and to restrict them with respect to a submanifold $G/H \
 (\hookrightarrow \widetilde{G}/\widetilde{H}$).

If $G/H$ is \lq\lq generic\rq\rq, namely,
 a {\bf {principal orbit}} in $\widetilde{G}/\widetilde{H}$
 in the sense of Richardson,
 then it is readily seen that
  discrete spectrum of the branching law $\pi|_{G}$
 gives a discrete series for $G/H$
 (\cite{$\text{Ko}_{10}$, \S 8};
  see also \cite{Hu, $\text{Ko}_{1,5}$, $\text{Li}_1$}
  for concrete examples).

However, some other interesting homogeneous spaces $G/H$ occur as
  non-principal orbits on $\widetilde{G}/\widetilde{H}$,
  where the above strategy does not work in general.
A remedy for this is to impose
  the {\bf admissibility of the restriction} of $\pi$,
  which justifies the restriction of $L^p$-functions to submanifolds,
   and then gives rise to many non-symmetric homogeneous
 spaces that admit discrete series representations.
 For example, let us consider the case where
$G =  \widetilde{G}^\tau$ and $H = \widetilde{G}^\sigma$
 for commuting involutive automorphisms  $\tau$ and $\sigma$
  of  $\widetilde{G}$
 such that $\widetilde{G}/\widetilde{H}$ satisfies (5.1).
 Then by using Theorem~C and an asymptotic estimate of invariant measures
  \cite{$\text{Ko}_6$},
  we have:
\proclaim{Theorem~J (discrete series for non-symmetric spaces,
\cite{$\text{Ko}_{10}$})}
Assume that there is $w \in W_{\sigma}$ such that
$$
   {\Bbb R}_+\operatorname{-span}\Delta^+(\frak p)_{\sigma,w}
    \cap \sqrt{-1} (\frak t_0^*)^{-\tau} = \{0\}.
\tag 5.2
$$
Then there exist infinitely many discrete series representations
 for any homogeneous space of $G$ that goes through $x \widetilde H \in
 \widetilde{G}/\widetilde{H}$ for any $x \in \widetilde{K}$.
\endproclaim
We refer to \cite{$\text{Ko}_{10}$, Theorem~5.1}
 for definitions of
 a finite group $W_\sigma$ and $\Delta^+(\frak p)_{\sigma, w}$.
The point here is that the condition (5.2) can be easily checked.

For instance,
 if $G \simeq Sp(2n, \Bbb R) \simeq \widetilde{G}/\widetilde{H}$
  (a group manifold),
 then Theorem~J implies that there exist discrete series %
 on all homogeneous spaces of the form:
$$
 G/H=Sp(2n,\Bbb R)/(Sp(n_0,\Bbb C) \times GL(n_1,\Bbb C) \times \dotsb \times
 GL(n_k,\Bbb C)),
 \quad
 (\sum n_i =n).
$$
The choice of $x$ in Theorem~J
 corresponds to the partition
  $(n_0, n_1, \dots, n_k)$.
We note that the above $G/H$ is a symmetric space
 if and only if
$
     n_1 = n_2 = \dots = n_k = 0.
$

The restriction of unitary representations gives new methods
 even for symmetric spaces
  where harmonic analysis has a long history of research.
Let us state two results that are proved by the theory
 of discretely decomposable restrictions.
\proclaim{Theorem~K (holomorphic discrete series for symmetric spaces)}
Suppose $G/H$ is a non-compact irreducible symmetric space.
Then {\rm (i)} and {\rm (ii)} are equivalent:
\item{{\rm i)}}
There exist unitary highest weight representations of $G$
 that can be realized as subrepresentations of $L^2(G/H)$.
\item{{\rm ii)}}
$G/K$ is Hermitian symmetric and
$H/(H \cap K)$ is its totally real submanifold.
\endproclaim

This theorem in the group manifold case is a restatement of
  Harish-Chandra's well-known result.
The implication (ii) $\Rightarrow$ (i)
 was previously obtained by a different geometric approach
 ('Olafsson-\O rsted \cite{Ol-\O}).
Our proof uses a general theory of discretely decomposable restrictions,
 especially, Theorems~B, H and J.

\proclaim{Theorem~L (exclusive law of discrete spectrum for
                restriction and induction)}
Let $G/G'$ be a non-compact irreducible symmetric space,
 and $\pi \in \widehat{G}$.
Then both (1) and (2) cannot occur simultaneously.
\item{{\rm 1)}}
The restriction $\pi|_{G'}$ is
 infinitesimally discretely decomposable.
\item{{\rm 2)}}
$\pi$ is a discrete series representation
  for the homogeneous space $G/G'$.
\endproclaim

We illustrate Theorems~K and L by $G=SL(2, \Bbb R)$.
The examples below are well-known results on harmonic analysis,
 however,
  the point is that they can be proved by a simple idea
   coming from restrictions of unitary representations.
\example{Example~5.4}
1) Holomorphic discrete series exist %
  for $G/H=SL(2, \Bbb R)/SO(1,1)$ (a hyperboloid of one sheet).
This is explained by Theorem~K
 because the geodesic $H/(H \cap K)$ is obviously totally real in
  the Poincar\'e disk $G/K = SL(2, \Bbb R)/SO(2)$.
\newline{2)}
There is no discrete series for
  the Poincar\'e disk $G/K = SL(2, \Bbb R)/SO(2)$.
This fact is explained by Theorem~L
  because any representation of $G$ is obviously
 discretely decomposable  when restricted to a compact $K$.
\endexample

\baselineskip 4.5mm \parindent 8mm

\specialhead \noindent \boldLARGE 6. Modular varieties, vanishing
theorem \endspecialhead

Retain the setting as in Definition~2.2.
Let
$
 \Gamma' \subset\ \Gamma
$
be cocompact torsion-free discrete subgroups
        of $G' \subset G$, respectively.
For simplicity,
 let $G'$ be a semisimple Lie group
 without compact factors.
Then,
 both of the double cosets
$
   X:=\Gamma \backslash G / K
 \text{ and }
   Y:=\Gamma' \backslash G' / K'
$
 are compact, orientable, locally Riemannian symmetric spaces.
Then,
 the inclusion $G' \hookrightarrow G$ induces a natural map
$
   \iota:Y \to X.
$
The image $\iota(Y)$ defines a totally geodesic submanifold in $X$.
Consider the induced homomorphism of the homology groups
 of degree $m:= \dim Y$,
$$
     \iota_*:H_m(Y;\Bbb Z)
               \rightarrow
       H_m(X;\Bbb Z).
$$
The {\bf {modular symbol}} is defined to be
 the image $\iota_*[Y] \in H_m(X;\Bbb Z)$ of
the fundamental class $[Y] \in H_m(Y;\Bbb Z)$.
Though its definition is simple,
 the understanding of modular symbols is highly non-trivial.

Let us first recall some results of
 Matsushima-Murakami and Borel-Wallach on the de Rham
 cohomology group $H^*(X; \Bbb C)$ summarized as:
$$
             H^* (X;\Bbb C)
                             =\bigoplus _{\pi \in \widehat{G}}
                              H^*(X)_\pi,
\quad
H^*(X)_\pi:= \operatorname{Hom}_G(\pi, L^2(\Gamma\backslash G)) \otimes
                 H^*(\frak g, K; \pi_K).
\tag 6.1
$$
The above result describes the topology of a single $X$
 by means of representation theory.
For the topology of the pair $(Y, X)$,
 we need restrictions of representations:
\proclaim{Theorem~M (vanishing theorem for modular symbols, \cite{Ko-Od})}
If
$$
     \operatorname{AS}_K(\pi) \cap \operatorname{Cone}(G')
     = \{0\}, \quad
     \pi \ne {\boldkey 1},
$$
then the modular symbol
 $\iota_* [Y] $
 is annihilated by the $\pi$-component $H^m(X)_\pi$
 in the perfect paring
$                 H^m (X; \Bbb C) \times
                  H_m(X ;\Bbb C) \to \Bbb C.
$
\endproclaim

Theorem~M determines, for example, the middle Hodge components
 of totally real modular symbols of compact Clifford-Klein forms
 of type IV domains.

The discreteness of irreducible decompositions plays
 a crucial role both in Matsushima-Murakami's formula (6.1)
 and in a vanishing theorem for modular varieties
 (Theorem~M).
In the former $L^2(\Gamma \backslash G)$ is $G$-admissible
 (Gelfand and Piateski-Shapiro),
 while the restriction $\pi|_{G'}$
 is $G'$-admissible
 (cf\. Theorem~B)
 in the latter.

\baselineskip 4.5mm \parindent 8mm

\specialhead \noindent \boldLARGE References \endspecialhead
\widestnumber\key{MMMMm}
\ref
    \key A
    \by M\. F\. Atiyah
    \paper The Harish-Chandra character
    \jour London Math\. Soc\. Lecture Note Series
    \vol 34
    \pages 176--181
    \yr 1979
\endref
\ref
    \key D-Vs
    \by M\. Duflo and J\. Vargas
    \jour in preparation
\endref
\ref
    \key FJ
    \by M\. Flensted-Jensen
    \paper Discrete series for semisimple symmetric spaces
    \jour Annals of Math\.
    \pages 253--311
    \vol 111
    \yr 1980
\endref
\ref
        \key Ge-Gv
        \by I\. M\. Gelfand and M\. I\. Graev
        \paper Geometry of homogeneous spaces, representations of groups
           in homogeneous spaces, and related questions of integral geometry
        \jour Transl\. II\. Ser., A\. M\. S\.
    \yr 1964
        \vol 37
        \pages   351--429
\endref
\ref
    \key $\text{Gr-W}_1$
    \by B\. Gross and N\.  Wallach
    \paper A distinguished family of unitary
           representations for the exceptional groups of real rank $= 4$
\jour Progress in Math\.
    \pages 289--304
    \publ Birkh\" auser
    \vol 123
    \yr 1994
\endref
\ref
\key $\text{Gr-W}_2$
\by B\. Gross and N\. Wallach
\paper Restriction of small discrete series representations
        to symmetric subgroups
\pages 255--272
\jour Proc\. Sympos\. Pure Math\.
\vol 68
\yr 2000
\publ A.M.S.
\endref
\ref
        \key $\text{Ho}_1$ %
        \by R\. Howe
        \paper $\theta$-series and invariant theory
        \jour Proc\. Sympos\. Pure Math\.
        \pages 275--285
        \publ A.M.S.
        \vol 33
        \yr 1979
\endref
\ref
   \key $\text{Ho}_2$ %
   \by R\. Howe
   \paper Wave front sets of representations of Lie groups
   \pages 117--140
   \jour Automorphic forms, representation theory, and arithmetic
\publ Tata
   \yr 1981
\endref
\ref
    \key Hu %
    \by J-S\. Huang
    \paper Harmonic analysis on compact polar homogeneous spaces
    \jour Pacific J\. Math\.
    \vol 175
    \yr 1996
    \pages 553---569
\endref
\ref
\key Hu-P-S
\by J-S. Huang, P. Pand\v zi\'c, and G. Savin
\paper New dual pair correspondences
\jour Duke Math\.
\vol 82
\yr 1996
\pages 447--471
\endref
\ref
   \key Hu-Vo
   \by J-S. Huang and D. Vogan
   \jour personal communications
   \yr 2001
\endref
\ref
    \key J-Vr
    \by H\. P\. Jakobsen and M\. Vergne
    \paper Restrictions and expansions of holomorphic representations
    \jour J\. Funct. Anal.
    \vol 34
    \yr 1979
    \pages 29--53
\endref
\ref
   \key Ks-Vr
   \by M\. Kashiwara and M\. Vergne
   \paper $K$-types and singular spectrum
    \inbook Lect\. Notes in Math\.
    \yr 1979, 177--200
   \publ Springer
    \vol 728
\endref
\ref
      \key Kn-Vo
      \by A\. Knapp and D\. Vogan, Jr\.
    \book Cohomological Induction and Unitary Representations
    \publ Princeton U\.P\.
    \yr 1995
\endref
\ref
        \key $\text{Ko}_1$
    \by T\. Kobayashi
        \paper Unitary representations realized in $L^2$-sections
              of vector bundles over semi-simple symmetric spaces
        \yr 1989
        \pages 39--54
        \jour Proc\. of the 27-28th Symp\. of Funct\.
          Anal\. and Real Anal\.
          \lang Japanese
    \publ Math\. Soc\. Japan
\endref
\ref
        \key $\text{Ko}_2$ %
       \by T\. Kobayashi
        \book Singular Unitary Representations and Discrete Series for
          Indefinite Stiefel Manifolds
        $U(p,q;{\Bbb F})/U(p-m,q;{\Bbb F})$
        \bookinfo Memoirs of A.M.S.
        \vol 462
        \yr 1992
\endref
     \ref
    \key $\text{Ko}_3$ %
  \by T\. Kobayashi
    \paper  The Restriction of $A_{\frak q}(\lambda)$ to reductive subgroups
    \jour Part I, Proc\. Japan Acad\.
    \yr 1993
    \vol 69
    \pages 262--267
    \moreref Part II, ibid. {\bf 71} 1995, 24--26
\endref
\ref
        \key $\text{Ko}_4$
    \by T\. Kobayashi
        \paper Harmonic analysis on homogeneous manifolds of reductive type
               and unitary representation theory
        \jour %
          Transl\., Series II,
         Selected Papers on Harmonic Analysis, Groups, and Invariants
\eds K. Nomizu
        \vol 183
        \yr  1998
        \pages 1--31
        \publ A.M.S.
\endref
\ref
     \key $\text{Ko}_5$ %
     \by T\. Kobayashi
     \paper  Discrete decomposability of the restriction
          of $A_{\frak q}(\lambda)$ with respect to reductive subgroups
          and its applications
        \jour Invent\. Math\.
        \yr 1994
    \vol 117
    \pages 181--205
\endref
\ref
    \key $\text{Ko}_6$ %
    \by T\. Kobayashi
   \paper  Invariant measures on  homogeneous manifolds of reductive type
   \jour J\. reine und angew\. Math\.
   \yr 1997
   \vol 490
   \pages 37--53
\endref
\ref
    \key $\text{Ko}_7$
  \by T\. Kobayashi
   \paper Discrete decomposability of the restriction of $A_\frak q(\lambda)$
          with respect to reductive subgroups {\rm II}
       ---  micro-local analysis and asymptotic $K$-support
    \jour Annals of Math.
    \yr 1998
    \pages 709--729
    \vol 147
\endref
\ref
    \key $\text{Ko}_8$
    \by T\. Kobayashi
    \paper Multiplicity free branching laws
            for unitary highest weight modules
    \yr 1997, 7--13
    \inbook Proceedings of the Symposium on Representation Theory held at Saga, Kyushu
    \ed K. Mimachi
    \endref
\ref
    \key $\text{Ko}_9$
  \by T\. Kobayashi
   \paper Discrete decomposability of the restriction of $A_\frak q(\lambda)$
          with respect to reductive subgroups {\rm III}
         --- restriction of Harish-Chandra modules and associated varieties
    \jour Invent\. Math\.
    \yr 1998
    \vol 131
    \pages 229--256
\endref
\ref
  \by T\. Kobayashi
   \key $\text{Ko}_{10}$
   \paper Discrete series representations for the orbit spaces
     arising from two involutions of real reductive Lie groups
    \jour J\. Funct\. Anal\.
    \vol 152
    \yr 1998
    \pages 100--135
\endref
\ref
   \key $\text{Ko}_{11}$
   \by T\. Kobayashi
   \paper    Discretely decomposable restrictions of
  unitary representations of reductive Lie groups
  ---examples and conjectures
   \yr 2000
   \pages 99-127
\jour Adv\. Stud\. Pure Math\.
   \vol 26
\endref
\ref
    \key $\text{Ko}_{12}$
    \by T\. Kobayashi
    \paper  Theory of discrete decomposable branching laws
         of unitary representations of semisimple Lie groups
         and some applications
    \jour Sugaku Exposition, Transl\. Ser\.
    \publ A.M.S.
    \toappear
\endref
\ref
    \key $\text{Ko}_{13}$
    \by T\. Kobayashi
    \jour in preparation
\endref
\ref
  \by T\. Kobayashi and T\. Oda
   \key Ko-Od
   \paper  Vanishing theorem of modular symbols on locally symmetric spaces
   \jour Comment\. Math\. Helvetici
   \yr 1998
   \vol 73
   \pages 45--70
\endref
\ref
   \key $\text{Ko-\O}_1$
   \by T\. Kobayashi and B\. \O rsted
    \paper Conformal geometry and branching laws for unitary representations
           attached to minimal nilpotent orbits
   \jour C\. R\. Acad\. Sci\. Paris
    \vol 326
    \pages 925--930
    \yr 1998
\endref
\ref
   \key $\text{Ko-\O}_2$
   \by T\. Kobayashi and B\. \O rsted
    \paper Analysis on the minimal representation of O(p,q),
     {\rm I, II, III}
   \jour preprint
\endref
\ref \key Le-Lo \by S-T. Lee and H-Y. Loke \paper Degenerate principal series of $U(p,q)$ and Spin$(p,q)$ \jour
preprint
\endref
\ref
    \key $\text{Li}_1$
    \by J-S\. Li
    \paper On the discrete series of generalized Stiefel manifolds
    \jour Trans\. A\.M\.S\.
    \yr 1993
    \pages 753--766
    \vol 340
\endref
\ref
\key   $\text{Li}_2$     \by J-S\. Li
    \paper Two reductive dual pairs in groups of type $E$
    \jour Manuscripta Math\.
    \vol 91
    \yr 1996
    \pages 163--177
\endref
\ref
\key $\text{Lo}_1$
\by H-Y\. Loke
\paper Restrictions of quaternionic representations
\jour J\. Funct\. Anal\.
\vol 172
\yr 2000
\pages 377--403
\endref
\ref
\key $\text{Lo}_2$
\by H-Y\. Loke
\paper Howe quotients of unitary characters and unitary lowest weight modules
\jour preprint
\endref
\ref
    \key Mg
    \by G\. Margulis
    \paper Existence of compact quotients of homogeneous spaces,
           measurably proper actions,
           and decay of matrix coefficients
    \jour Bul\. Soc\. Math\. France
    \vol 125
    \yr 1997
    \pages 1--10
\endref
\ref
    \key Mt
    \by S\. Martens
    \paper The characters of the holomorphic discrete series
    \jour Proc\. Nat\. Acad\. Sci\. USA
    \vol 72
    \yr 1975
    \pages 3275-3276
\endref
\ref
    \key Mk-Os
    \by T\. Matsuki and T\. Oshima
    \paper A description of discrete series for semisimple symmetric spaces
    \jour Adv\. Stud\. Pure Math\.
    \pages 331--390
    \vol 4
    \yr 1984
\endref
\ref
\key N-Oc-T
\by K\. Nishiyama, H\. Ochiai, and K\. Taniguchi
\paper Bernstein degree and associated cycles of Harish-Chandra modules
  --- Hermitian symmetric case
\jour Asterisque
\vol 273
\yr 2001
\pages 13--80
\endref
\ref
     \key Oh
    \by H\. Oh
    \paper Tempered subgroups and representations with
           minimal decay of matrix coefficients
    \jour Bull\. Soc\. Math\. France
    \vol 126
    \yr 1998
    \pages  355--380
\endref
\ref
    \key Ol-\O
    \by G\. \'Olafsson and B\. \O rsted
   \paper The holomorphic discrete series of an affine symmetric space, {\rm I}     \jour J\. Funct\. Anal\.
    \yr 1988
    \pages 126--159
    \vol 81
\endref
\ref
    \key \O-Vs
    \by B\. \O rsted and J\. Vargas
    \paper Restriction of square integrable representations:
           discrete spectrum
    \jour preprint
\endref
\ref
    \key Vo
    \by D\. Vogan, Jr\.
    \paper Associated varieties and unipotent representations
        \pages 315--388
\jour Progress in Math\.
    \publ Birkh\" auser
    \vol 101
    \yr 1991
\endref
\ref
    \key Vo-Z
    \by D\. Vogan, Jr\. and G\. Zuckerman
    \paper Unitary representations with non-zero cohomology
    \jour Compositio Math\.
    \vol 53
    \pages 51--90
    \yr 1984
\endref
\ref
    \key X
    \by J\. Xie
    \paper Restriction of discrete series of $SU(2,1)$ to
  $S(U(1) \times U(1,1))$
    \jour J\. Funct\. Anal\.
    \yr 1994
    \vol 122
    \pages 478--518
    \finalinfo ph.D\. dissertation, Rutgers University %
\endref
\enddocument